\documentclass[preprint]{elsarticle}

\usepackage{lineno,hyperref}
\modulolinenumbers[5]

\usepackage{latexsym,amsfonts,amssymb,amsmath,amsthm}
\usepackage{graphicx}
\usepackage{braket}
\usepackage{bbm}
\usepackage{amsmath}

\usepackage[usenames,dvipsnames]{color}
\usepackage{ulem}

\begin{document}


\newtheorem{theo}{\textbf{\ \ \quad Theorem}}[section]
\newtheorem{li}{\textbf{\ \ \quad Li}}[section]
\newtheorem{lem}{\textbf{\ \ \quad Lemma}}[section]
\newtheorem{remark}{\textbf{\ \ \quad Remark}}[section]
\newtheorem{col}{\textbf{\ \ \quad Corollary}}[section]
\newtheorem{prop}{\textbf{\ \ \quad Proposition}}[section]
\newtheorem{defi}{\textbf{\ \ \quad Definition}}[section]

\newcommand{\lbl}[1]{\label{#1}}
\newcommand{\bib}[1]{\bibitem{#1}\qquad\framebox{\scriptsize #1}}
\newcommand{\bth}[2]{\vskip 8pt\noindent\bf #1\hskip 2pt\bf#2\it \hskip 8pt}
\newcommand{\be}{\begin{equation}}
\newcommand{\ee}{\end{equation}}
\newcommand\bes{\begin{eqnarray}}
\newcommand\ees{\end{eqnarray}}
\newcommand{\bess}{\begin{eqnarray*}}
\newcommand{\eess}{\end{eqnarray*}}
\newcommand\bee{{\mbox{\boldmath $e$}}}
\newcommand\bR{{\mathbf R}}
\newcommand\bC{{\mathbf C}}\newcommand{\rhp}{\rightharpoonup}
\newcommand\RE{{\rm Re}}\newcommand\mA{{\rm \bf (A)}}
\newcommand\mB{{\rm \bf (B)}}
\newcommand{\nm}{\nonumber}
\newcommand{\p}{\partial}
 \newcommand{\dd}{\displaystyle}
 \newcommand{\rd}{{\rm d}}
\newcommand{\dis}{\mbox{\rm d}_\ep\mbox{}}
\newcommand\ep{{\varepsilon}}
\newcommand\mean[1]{{\langle #1\rangle}}
\newcommand\meanep[1]{{\langle #1\rangle_\ep}}
\newcommand\pd[2]{{\frac{\partial #1}{\partial #2}}}
\newcommand\cc{{\mbox{\boldmath $c$}}}
\newcommand\Mani{{\mathcal M}}
\newcommand\NM{{N_\xi{\mathcal M}}}
\newcommand\lm{\lambda}
\newcommand{\ol}{\overline}
\newcommand{\ds}{\displaystyle}
\newcommand{\ul}{\underline}
\newcommand{\na}{\nabla}
\newcommand{\vmm}{\vspace{1.5mm}}
\newcommand\bI{{\mathbf I}}
\newcommand\bII{{\mathbf {II}}}
\renewcommand\dim{{\hbox{\rm dim}}}
\newcommand\sgn{{\hbox{\rm sgn}}}
\newcommand\trace{{\hbox{Tr}}}
\newcommand\defequal{{ \ \stackrel{\rm def}{=} \ }}
\newcommand\bu{{\mathbf u}}
\newcommand\bv{{\mathbf v}}
\newcommand\bc{{\mathbf c}}
\newcommand\bF{{\mathbf F}}
\newcommand\bG{{\mathbf G}}
\newcommand\bA{{\mathbf A}}
\newcommand\bB{{\mathbf B}}
\newcommand\bU{{\mathbf U}}
\newcommand\bY{{\mathbf Y}}
\newcommand\bX{{\mathbf X}}
\newcommand\bD{{\mathbf D}}
\newcommand\bK{{\mathbf K}}
\newcommand\R{{\mathbb R}}
\newcommand\F{{\mathcal F}}
\newcommand\K{{\mathcal K}}

\begin{frontmatter}

\title{A regularity result  for the fractional Fokker-Planck equation with Ornstein-Uhlenbeck drift }

\author[mymainaddress]{Xiaoxia Xie\corref{mycorrespondingauthor}}
\cortext[mycorrespondingauthor]{Corresponding author}
\ead{xxie12@iit.edu}

\author[mymainaddress]{Jinqiao Duan}
\ead{duan@iit.edu}

\author[mymainaddress]{Xiaofan Li}
\ead{lix@iit.edu}

\author[mysecondaryaddress]{Guangying Lv}
\ead{gylvmaths@henu.edu.cn}

\address[mymainaddress]{Department of Mathematics, Illinois Institute of Technology, Chicago, IL 60616, United States}
\address[mysecondaryaddress]{ Institute of Contemporary Mathematics, Henan University, Kaifeng, Henan 475001,  China}

\begin{abstract}
Despite there are numerous theoretical studies of stochastic differential equations with a symmetric $\alpha$-stable L\'evy  noise, very few regularity results exist in the case of $0<\alpha\leq1$. In this paper, we study the fractional Fokker-Planck equation 
with  Ornstein-Uhlenbeck drift, and prove that there exists a unique solution, which is $C^\infty$ in space for $t>0$ when $\alpha\in (0, 2]$.
\end{abstract}


\end{frontmatter}


\section{Introduction}

In this paper, we are interested in the regularity of solutions of the generalized Fokker-Plank equation  for $ x\in \R^d$ and $t>0$
\be
\label{V-drift-F-P-alpha}
\p_t u=\mathcal I[u]-\nabla\cdot( f(x)u), \quad u(0, x)=u_0(x),
\ee
where $\mathcal I$ is an operator defined in \eqref{I},  $f(x):\R^d\to \R^d$ in the drift term is a deterministic driving force, and $u_0$ in $ L^\infty (\R^d)$ is nonnegative. We assume without loss of generality that $\int_{\R^d}u_0(x)dx=1$.   When $\mathcal I$ is the classical Laplacian, Eq \eqref{V-drift-F-P-alpha} is the classical Fokker-Plank equation. If $f(\cdot)$ fulfills the usual assumptions of regularity needed in the Laplacian perturbation theory, it is well-known that the solution is smooth  due to the  regularizing effect of the parabolic operator. If $f(\cdot)$ is less regular, there is a series of papers by Lions, {\it et al.} devoted to vector field $f(\cdot)$ having bounded divergence and some Sobolev type regularity, and it is shown that the $H^1$ regularity in $x$ is sharp for $f(x)$ being in Sobolev spaces instead of being Lipschitz or continuously bounded (see \cite{B-L,B-L-renorm, D-L, D-L-cauchy}). 

The Laplacian is the infinitesimal generator of the Brownian motion, hence the classical Fokker-Planck equation   corresponds to a stochastic system with Gaussian noises. However, the noise is often non-Gaussian, and in this paper, we focus on the symmetric $\alpha$-stable L\'evy motion, for which the infinitesimal generator is the fractional Laplacian, denoted by  $\Delta^{\alpha/2}$ with $\alpha\in(0, 2)$. Eq \eqref{V-drift-F-P-alpha}  with $\mathcal I=\Delta^{\alpha/2}$ is called the fractional Fokker-Plank equation (see Section 2 for details).

There are many results (see \cite{B-J}, \cite{B-S-S}-\cite{D-I}, \cite{ G-D-L}, \cite{ H-D-G}-\cite{K-N},  \cite{SL}-\cite{W-T},\cite{S} and references therein) on gradient perturbation of the fractional Laplacian operator in the frame work of regular vector field, i.e.  the well posedness and regularity of the solution to the equation $u_t=\Delta^{\alpha/2}u-f\cdot \nabla u$ with $u_0\in L^\infty(\R^d)$ and $f$ satisfying some regularity assumptions. Depending on the value of $\alpha$, there are three cases. {First},  $\alpha\in (0, 1)$ is known as the supercritical case because $\Delta^{\alpha/2}u$ is of lower order than $f\cdot \nabla u$. In  \cite{SL}, it is shown that  the solution becomes differentiable with H\"older continuous derivatives, if $f\in C^{1-\alpha+\beta}$ for any $\beta\in (0,\alpha $). With a slightly weaker regularity assumption on $f$, it is shown in \cite{SL-holder} that $u$  is H\"older continuous if $f\in C^{1-\alpha}$ in space. {Second}, in the critical case of $\alpha=1$,  Silvestre in \cite{SL-holder} shows that $u$ becomes H\"older continuous for positive time  if $f\in L^{\infty}(\R^d)$, while the authors in \cite{C-S-regular} and \cite{K-N} achieved the same result by different approaches 
when $f$ belongs to the BMO class and $\text{div}f=0$.
{Third}, when $\alpha\in (1, 2)$, the term $\Delta^{\alpha/2}u$ is  dominant, which is called the subcritical case,  and  better regularities  have been proved under more general conditions. Actually, in the one-dimensional case ($d=1$), Droniou, {\it et al.} in \cite{D-G-V} show  that the unique solution exists and $u\in C^\infty((0, \infty), \R)$  for $f\cdot \nabla u$ replaced by $\p_x(f(u))$. Their result is generalized to multidimensional cases  in \cite{I} with the term $H(t, x, u, \nabla u)$ instead of $\p_x(f(u))$. Note that $f(u)$ and $H(t, x, u, \nabla u)$ must satisfy certain conditions in order for $\Delta^{\alpha/2}$ to be the dominate term. The estimate of the fractional heat kernel  perturbed by gradient operator is obtained  in \cite{B-J} if  $f$ belongs to a Kato class. See also \cite{Chen-K, B-S-S, B-G} for two-sided fractional heat kernel estimates.

However, for the vector field $f$ lying in some Sobolev spaces instead of being uniformly bounded in $x$,  the drift term could be stronger than the diffusion term $\mathcal I[u]$ even in the case of  $\mathcal I=\Delta$. Hence  it would not be possible to prove any regularity for $u$ in general and there are relatively fewer papers. For the fractional Fokker-Plank equation, the existence of weak $L^p$-solution was recently obtained in \cite{W-T}, if $f(\cdot)\in L_{loc}^q(\R^d)$, where $q$ is the conjugate exponent of $p$. The solution is also unique with further assumptions on $f(\cdot)$. However, there is no regularity result so far.

In the rest of the paper, we focus on the fractional Fokker-Plank equation with the Ornstein-Uhlenbeck drift, i.e., Eq \eqref{V-drift-F-P-alpha} with $\mathcal I=\Delta^{\alpha/2}$ and $f(x)=-x$,
\be
\label{OU-drift-F-P-alpha}
\p_t u=\Delta^{\alpha/2}u-\nabla\cdot(- xu), \quad 
u(0, x)=u_0(x).
\ee
In this case, the unique  weak $L^p$-solution exists by applying the results in \cite{W-T}. However, we show that the unique solution to Eq \eqref{OU-drift-F-P-alpha} is smooth in this work. More precisely, the solution of \eqref{OU-drift-F-P-alpha} can be represented as
 \be
\label{V-drift-F-P-alpha-sol}
u(t, x):=( p\ast u_0)(x)=\int_{\R^d}p(t, x-y)u_0(y)dy,
\ee
where $p(t, x, y)=p(t, y-x)$ is the fundamental solution of \eqref{OU-drift-F-P-alpha} (or $p$ is  called the heat kernel of the fractional Laplacian perturbed by the  Ornstein and Uhlenbeck drift). One can easily check that 
 \begin{equation}
\label{p-hatp=}
 p(t, x, y)=e^{dt}\widehat p\left(\frac{e^{\alpha t}-1}{\alpha}, e^{t}x, y\right),
\end{equation}
where $\widehat p(t, x, y)=\widehat p(t, y-x)$ is fundamental solution to the fractional heat equation, i.e.,  $\p_t u=\Delta^{\alpha/2}u$ (or  $\widehat p$ is called the fractional heat kernel).
 It is well-known that $\widehat p(t, \cdot)$  is $C^\infty$  for all $t>0$  when $\alpha\in (0, 2)$ (see \cite{A}).
Due to \eqref{p-hatp=}, the solution of \eqref{OU-drift-F-P-alpha} has the same regularity as that the fractional heat equation, as shown in the following theorem.
 \begin{theo}
\label{thm-OU-drift-u-infty}
  Let $u_0\in L^\infty(\R^d)$ and be continuous almost everywhere. The function $u$ defined in \eqref{V-drift-F-P-alpha-sol} is the solution  to \eqref{OU-drift-F-P-alpha}. Furthermore, for $T>0$ or $T=\infty$, 
\be
\label{OU-drift-u-infty}
u(t, \cdot)\in C^\infty(\R^d), \quad \forall t\in (0, T].
\ee
\end{theo}
We remark that from the stochastic view point, the $\alpha$-stable  Ornstein-Uhlenbeck  process was first  studied in the case of $\alpha=1$ in \cite{G-O}, and the sharp estimates of the mean first exist time from the ball was given for $0<\alpha<2$  in \cite{J-est-mean}, and it is shown  in \cite{J-harnack} that the Harnack inequality holds if $1\leq \alpha<2$ or $\alpha<1=d$. From the pseudo differential equations point of view, see \cite{G-I} for the exponential decay in time of the solution to the steady state of the associated stationary equation, and \cite{G-I2} for the applications of the L\'evy logarithmic Sobolev inequality to the study of the regularity of the solutions of the fractional heat equation and the asymptotic behavior of the L\'evy-Ornstein-Uhlenbeck process. For numerical results, see \cite{G-D-L}.  Finally, we refer to \cite{G-O, J-M-F, W, W-G-M-N} for applications of the fractional Ornstein-Uhlenbeck  process  in physical sciences.

We point out that our result in Theorem \ref{thm-OU-drift-u-infty} is a little surprising, especially for the supercritical case ($\alpha\in (0, 1)$), because it implies that  the fractional Laplacian operator $\Delta^{\alpha/2}$ has the smoothing effect like the Laplacian in the Ornstein-Uhlenbeck  drift case. We wish to expand the result to include more general drift terms in further studies.
The paper is organized as follows: in Section 2,  we  introduce some notations and definitions. And the proof of Theorem \ref{thm-OU-drift-u-infty} will be given in Section 3.
\section{Preliminaries}
 \setcounter{equation}{0}

Let us briefly recall  basic definitions of the {isotropic} $\alpha$-stable L\'evy operator and the potential theory of the fractional Laplacian (see\cite{A, B-G, B-H, L}).

The L\'evy-Khinchine formula implies that there exists a symmetric semi-definite $d\times d$ constant matrix $\sigma=(\sigma_{i, j})$, a constant vector $b=(b_i)\in \R^d$ and a nonnegative measure $\nu$ on $\R^d$ such that the L\'evy operator writes:
\be
\label{I}
\mathcal I[u](x)=\Braket{\sigma \nabla u; \nabla u}+\Braket{b;\nabla u}+\int_{\R^d} (u(x+z)-u(x)-\nabla u(x)\cdot z {\mathbbm 1}_{\{|z|<1\}})\nu(dz),
\ee
where $\Braket{\cdot;\cdot}$ is the inner product.
 The so-called L\'evy measure $\nu$ can be singular at the origin and satisfies $\displaystyle\nu(\{0\})=0$  and $\displaystyle \int\min(1, |z|^2)\nu(dz)<\infty$
 and the matrix $\sigma$ characterizes the diffusion (or Gaussian) part of the operator (it can be null if $\sigma=0$), while $b$ and $\nu$ characterizes the drift part and the pure jump part, respectively. The characteristic exponent of the L\'evy process is denoted by
$\displaystyle \psi(\xi)=\Braket{\sigma \xi; \xi}+\Braket{b; \xi}+a(\xi)$,
  where $\displaystyle a(\xi)=\int_{\R^d}[1-\cos(\xi\cdot z)]\nu(dz)$ is associated with the pure jump part of the underlying process.
 The symmetric $\alpha$-stable L\'evy process refers to the case
$b=0$,  $\sigma=0$,  and  $\displaystyle \nu(dz)=\mathcal A_{d, \alpha}\frac{dz}{|z|^{d+\alpha}}$,
where $\mathcal A_{d, \alpha}$ is so chosen that $a(\xi)=|\xi|^\alpha$  for $\xi \in \R^d$ (see \cite{B-J} and references therein).
 
 For $\alpha\in (0, 2)$ and $\psi\in C_c^\infty(\R^d)$, the fractional Laplacian is
\be
\label{delta-alpha-2}
\Delta^{\frac{\alpha}2}\phi(x):=A_{d,\alpha,} \lim_{\epsilon\to 0} \int_{|z|>\epsilon}\frac{\phi(x+z)-\phi(x)}{|z|^{d+\alpha}}dz.
\ee
Note that if $\alpha=2$ the L\'evy operator has no jump part, and  $\mathcal I[\cdot]$ is the Laplacian operator in this case. 

We let $\widehat p(t, \cdot)$ be the smooth real-valued function on $\R^d$ with Fourier transform
\be
\label{pro-tran-den-hatp}
\int_{\R^d}\widehat p(t, z)e^{iz\cdot \xi}dz=e^{-t|\xi|^\alpha}, \quad  t>0, \,\ \xi \in \R^d.
\ee
It follows from \eqref{pro-tran-den-hatp} that $\widehat p$ is the heat kernel of the fractional Laplacian and it is well-known that  for $t>0$, $x, y\in \R^d$, $x\neq y$, the sharp estimate of $\widehat p(t, x)$ is
\be
\label{Est-alpha}
\widehat p(t, x, y)\approx \min\left( \frac{t}{|x-y|^{d+\alpha}}, t^{-d/\alpha}\right)
\ee
(see \cite{B-J, B-S-S, D-G-V, I}),  where the notation $\approx $ means that either ratio of the sides is bounded by $C\in (0, \infty)$, and $C$ does not depend on the variables shown, here $t$ and $x$. 
The  estimate of the first order derivative of $\widehat p(t, x)$ is  was derived in  \cite[Lemma 5]{B-J},
 and we extend it to get the estimate of the $m$-th order derivative of $\widehat p(t, x)$ for any $m\geq 0$ as follows
\be
\label{Est-alpha-m=}
\p_x^m\widehat p(t, x)\approx\sum_{n=0}^{n=\lfloor \frac{m}{2}\rfloor}C_n |x-y|^{m-2n} \min \left\{ \frac{t}{|x-y|^{d+\alpha+2(m-n)}}, t^{-\frac{d+2(m-n)}{\alpha}}\right\},
\ee
where $\lfloor \frac m 2\rfloor$ means the largest integer that is less  than $\frac m2$.

\section{Proof of Main Result}
\setcounter{equation}{0}
\begin{proof}
To prove  Theorem  \ref{thm-OU-drift-u-infty}, we first prove the following result:  for $u_0\in L^\infty(\R^d)$ and $u_0(x)$ be continuous almost everywhere. The function $u(t, x)$ defined as
\be
\label{Free-F-P-alpha-sol}
u(t, x)=\int_{\R^d}\widehat p(t, y-x)u_0(y)dy
\ee
is the solution to $u_t=\Delta^{\alpha/2}u$ with $u(0, x)=u_0(x)$, and for $T>0$ or $T=\infty$,  we have
\be
\label{free-u-infty}
u(t, \cdot)\in C^\infty(\R^d), \quad \forall t\in (0, T].
\ee
 Firstly,  we show that $u(t, \cdot)\in C^\infty(\R^d)$ for all $t\in (0, T]$. Since $u_0\in L^\infty(\R^d)$, and  $\widehat p(t, \cdot)\in L^1(\R^d)$ for $t>0$, $\widehat p(t,\cdot)\ast u_0$ is well-defined. By the Young's inequality for the convolution in \eqref{Free-F-P-alpha-sol}, we have
\be
\label{free-sol-0}
\forall (t, x)\in (0, \infty)\times \R^d, \quad \|\widehat p(t, \cdot)\ast u_0(x)\|_{L^1(\R^d)}\leq \|u_0\|_{L^\infty(\R^d)}.
\ee
By the smoothness of $\widehat p(t, x, y)$, to show \eqref{free-u-infty}, it is sufficient to show that
\[
|\nabla_x^m\widehat p(t, x-y)|\in L^1(\R^d), \quad \forall m\geq 0, t_0<t<T, \text{ and } x, y\in \R^d.
\]
 Indeed, by \eqref{Est-alpha-m=},  we have
\begin{align*}
\int_{\R^d}|\nabla_x^m\widehat p(t, x-y)|dy
&\leq\sum_{n=0}^{n=\lfloor \frac m2 \rfloor} \int_{\R^d}|x-y|^{m-2n} \min\left\{ \frac{ t}{|x-y|^{d+\alpha+2(m-n)}}, { t^{-\frac{d+2(m-n)}{\alpha}}}\right\}dy\\
&=\sum_{n=0}^{n=\lfloor \frac m2 \rfloor}\int_{B(x, r)} |x-y|^{m-2n}{ t^{-\frac{d+2(m-n)}{\alpha}}}dy +\sum_{n=0}^{n=\lfloor \frac m2 \rfloor}\int_{B(x, r)^c}\frac{ t}{|x-y|^{d+\alpha+m}}dy\\
&\leq \sum_{n=0}^{n=\lfloor \frac m2 \rfloor} { t_0^{-\frac{d+2(m-n)}{\alpha}}}\int_{B(x, r)} |x-y|^{m-2n}dy+\sum_{n=0}^{n=\lfloor \frac m2 \rfloor}{ T}\int_{B(x, r)^c}\frac{1}{|x-y|^{d+\alpha+m}}dy\\
&<\infty.
\end{align*}
Sine $u_0\in L^\infty(\R^d)$ and $\widehat p(t, \cdot)\in C^{\infty}(\R^d)$, the theorem of continuity under the integral sign gives \eqref{free-u-infty}. The uniqueness of the solution $u(t, x)$ defined in \eqref{Free-F-P-alpha-sol} follows from the Fourier Transform.
And it is clear that for $x^0$ at which $u_0(x)$ is continuous, we have
\be
\label{Free-F-P-sol-conti-initial}
\lim_{t\to 0, x\to x_0}u(t, x)=u_0(x^0).
\ee
Indeed, fix $x^0\in\R^d, \epsilon>0$. { Suppose $u_0(x)$ is continuous at $x_0$}. Choose $\delta>0$ such that
\[
|u_0(y)-u(x^0)|<\epsilon \quad \text{ if } |y-x^0|<\delta, \text{ and } y\in \R^d.
\]
Then if $|x-x^0|<\frac{\delta}2$, we have
\begin{align*}
|u(t, x)-u_0(x^0)|\leq  &\int_{\R^d}\widehat p(t, x-y)|u_0(y)-u_0(x^0)|dy\\
= &\int_{B(x^0, \delta)}\widehat p (t, x-y)|u_0(y)-u_0(x^0)|dy+\int_{B(x^0, \delta)^c}\widehat p(t, x- y)|u_0(y)-u_0(x^0)|dy\\
\le &\epsilon+2\|u_0\|_{L^\infty}\int_{B(x^0, \delta)^c} \frac{t}{|x-y|^{d+\alpha}}dy.
\end{align*}
In the second term, we have  $|y-x^0|>\delta$. Hence
\[|y-x^0|\leq |y-x|+|x^0-x|<|y-x|+\frac{\delta}2\leq |y-x|+\frac12|y-x^0|.
\]
Thus $ |y-x|\geq\frac12 |y-x^0|$. Consequently, the second term in the last line tends to $0$ as $t\to 0^+$.
Hence, if $|x-x^0|<\frac{\delta}{2}$, and $t>0$ is small enough, $|u(x, t)-u_0(x^0)|<2\epsilon$.

Second, set $\tilde t=\frac{e^{\alpha t}-1}{\alpha}$ and $\tilde x=e^{ t} x$. By \eqref{p-hatp=}, we have
 \be
\label{p=hatp-nabla}
\nabla_x^m p(t, x, y)=e^{dt}e^{m t}\nabla_{\tilde x}^m\widehat p(\tilde t, \tilde x, y),
\ee
Hence by \eqref{Est-alpha-m=},  we have
\begin{align*}
\int_{\R^d}|\nabla_x^m p(t, x-y)|dy&\leq e^{(d+m)t}\sum_{n=0}^{n=\lfloor \frac m2 \rfloor} \int_{\R^d}|\tilde x-y|^{m-2n} \min\left\{ \frac{ \tilde t}{|\tilde x-y|^{d+\alpha+2(m-n)}}, { \tilde t^{-\frac{d+2(m-n)}{\alpha}}}\right\}dy\\
&=e^{(d+m)t}\left[\sum_{n=0}^{n=\lfloor \frac m2 \rfloor}\int_{B(\tilde x, r)} |\tilde x-y|^{m-2n}{ \tilde t^{-\frac{d+2(m-n)}{\alpha}}}dy +\sum_{n=0}^{n=\lfloor \frac m2 \rfloor}\int_{B( \tilde x, r)^c}\frac{\tilde t}{| \tilde x-y|^{d+\alpha+m}}dy\right]\\
&\leq e^{(d+m)t} \sum_{n=0}^{n=\lfloor \frac m2 \rfloor} { \tilde t^{-\frac{d+2(m-n)}{\alpha}}}\int_{B(\tilde x, r)} |\tilde x-y|^{m-2n}dy+\sum_{n=0}^{n=\lfloor \frac m2 \rfloor}{ \tilde t}\int_{B(\tilde x, r)^c}\frac{1}{| \tilde x-y|^{d+\alpha+m}}dy\\
&<\infty
\end{align*}
for any $0<t_0<t<T$. Hence, the theorem follows.
\end{proof}

\section*{Acknowledgement}

This work was partly supported by ANL Grant 3J-30001-0025A.


\section*{References}











\begin{thebibliography}{99}

\bibitem{A} D.Applebaum,  L\'evy processes and stochastic calculus, Cambridge university press, 2009.



%
%



\bibitem{B-G} R.M. Blumenthal,  R.K. Getoor,
Some theorems on stable processes,
Trans. Amer. Math. Soc. 95 1960 263-273. 


\bibitem{B-H} J.Bliedtner, W.Hansen,  Potential theory: an analytic and probabilistic approach to balayage, Springer Science $\&$ Business Media, 2012.
%
%
\bibitem{B-J}  K.Bogdan, T.Jakubowski, 
Estimates of heat kernel of fractional Laplacian perturbed by gradient operators, Comm. Math. Phys. 271 (2007), no. 1, 179-198.

\bibitem{B-L}C. Le Bris, P.-L. Lions,
Existence and uniqueness of solutions to Fokker-Planck type equations with irregular coefficients,
Comm. Partial Differential Equations 33 (2008), no. 7-9, 1272-1317.
\bibitem{B-L-renorm}C. Le Bris, P.-L. Lions, Renormalized solutions of some transport equations with partially $W^{1, 1}$ velocities and applications,  Annali di Matematica pura ed applicata 183, no. 1 (2004),  97-130.

\bibitem{B-S-S} K. Bogdan, A. St$\acute{o}$s, P. Sztonyk, 
Harnack inequality for stable processes on $d$-sets,
Studia Math. 158 (2003), no. 2, 163-198.





\bibitem{C-S-regular}L. Caffarelli,  L.  Silvestre,  Regularity theory for fully nonlinear integro-differential equations,  Communications on Pure and Applied Mathematics  62, no. 5 (2009) 597-638.
%
\bibitem{Chen-K} Z.-Q. Chen, T. Kumagai, Heat kernel estimates for stable-like processes on $d$-sets,  Stoch.
Process. Appl. 108 (2003), no. 1, 27-62.




\bibitem{D-G-V} J. Droniou, T. Gallou$\ddot{\text{e}}$t, J. Vovelle, 
Global solution and smoothing effect for a non-local regularization of a hyperbolic equation,
J. Evol. Equ. 3 (2003), no. 3, 499-521.
\bibitem{D-I}J. Droniou, C. Imbert, 
Fractal first-order partial differential equations,
Arch. Ration. Mech. Anal. 182 (2006), no. 2, 299-331.

\bibitem{D-L}R. J. DiPerna, P.-L. Lions,
Ordinary differential equations, transport theory and Sobolev spaces,
Invent. Math. 98 (1989), no. 3, 511-547. 

\bibitem{D-L-cauchy}R.J. DiPerna, P.-L. Lions, 
On the Cauchy problem for Boltzmann equations: global existence and weak stability, 
Ann. of Math. (2) 130 (1989), no. 2, 321-366. 





\bibitem{G-D-L}T. Gao, J. Duan, X. Li,  Fokker-Plank equations for stochastic dynamical systems with symmetric L$\acute{\text e}$vy motions, submitted.

\bibitem{G-I} I. Gentil, C. Imbert,  The L$\acute{\text e}$vy-Fokker-Planck equation: $\Phi$-entropies and convergence to equilibrium, Asymptot. Anal. 59 (2008), no. 3-4, 125-138.

\bibitem{G-I2} I. Gentil, C. Imbert,  Logarithmic Sobolev inequalities: regularizing effect of L$\acute{\text e}$vy operators and asymptotic convergence in the L$\acute{\text e}$vy-Fokker-Planck equation, Stochastics: An International Journal of Probability and Stochastics Processes 81, no. 3-4 (2009) 401-414.



\bibitem{G-O} P. Garbaczewski, R.  Olkiewicz,  Ornstein-Uhlenbeck-Cauchy process,
J. Math. Phys. 41 (2000), no. 10, 6843-6860.


\bibitem{H-D-G} J. He, J. Duan, H. Gao,  A nonlocal Fokker-Plank equation for non-Gaussian  stochastic dynamical systems, To appear in {\it Applied Math. Lett.}, 2015.
    
    
\bibitem{I} C. Imbert, A non-local regularization of first order Hamilton-Jacobi equations, J. Differential Equations 211 (2005), no. 1, 218-246.



\bibitem{J-est-mean}T. Jakubowski, 
The estimates of the mean first exit time from a ball for the $\alpha$-stable Ornstein-Uhlenbeck processes. Stochastic Process, Appl. 117 (2007), no. 10, 1540-1560.

\bibitem{J-harnack} T. Jakubowski,  On Harnack inequality for $\alpha$-stable Ornstein-Uhlenbeck processes, Math. Z. 258 (2008), no. 3, 609-628.

\bibitem{K-N} A. Kiselev, F. L. Nazarov,  Variation on a theme of Caffarelli and Vasseur, Journal of Mathematical Sciences 166, no. 1 (2010) 31-39.


\bibitem{L} N. M. Landkof,  Foundations of modern potential theory. Vol. 180. Springer, 1972.
\bibitem{SL}L. Silvestre, On the differentiability of the solution to an equation with drift and fractional diffusion,
Indiana University Mathematical Journal. 61 (2012), no. 2, 557-584.

\bibitem{SL-holder}L. Silvestre,  H\"older estimates for advection fractional-diffusion equations,
Annali della Scuola Normale Superiore di Pisa. Classe di Scienze. Accepted for publication.
\bibitem{W-T} J. Wei, R. Tian, Well-posedness for the fractional Fokker-Plank equations, submitted.


\bibitem{J-M-F} S. Jespersen, R.  Metzler, H.C. Fogedby,   L$\acute{\text e}$vy flights in external force fields: Langevin and fractional Fokker-Planck equations and their solutions,
Phys. Rev. E 59, 2736.



\bibitem{S} A. St$\acute{\text o}$s, Symmetric $\alpha$-stable processes on $d$-sets, Bull. Polish Acad. Sci. Math. 48 (2000), no. 3, 237-245.

%




%
\bibitem{W}W. A. Woyczy\'nski, 
L$\acute{\text e}$vy processes in the physical sciences, L$\acute{\text e}$vy processes, 241-266, Birkh$\ddot{\text a}$user Boston, Boston, MA, 2001.

\bibitem{W-G-M-N} B.J. West, P. Grigolini, R. Metzler, T. F. Nonnenmacher, 
Fractional diffusion and L$\acute{\text e}$vy stable processes,
Phys. Rev. E (3) 55 (1997), no. 1, part A, 99–106.



\end{thebibliography}
 \end{document}